\documentclass[11pt,leqno]{article}

\usepackage{amsmath} 
\usepackage{amsfonts} 
\usepackage{amssymb} 
\newtheorem{theo}{Theorem} 
\newtheorem{cor}{Corollary} 
\newtheorem{lemma}{Lemma} 
\usepackage[hyperfootnotes=false,colorlinks]{hyperref}

\begin{document}

\title{Quantiles as minimizers} 
\date{November 2014} 
\author{Michel {\sc Valadier}\thanks 
{\ Retired Professor, Math\'ematiques, Universit\'e Montpellier II, 
place Eug\`ene Bataillon, Case courier 051, 34095 Monptellier Cedex, France. 
email: mivaladier{\string @}wanadoo.fr}} 

\maketitle 

\begin{abstract} 
A real random variable admits median(s) and quantiles. 
These values minimize convex functions on $\mathbb R$. 
We show by ``Convex Analysis'' arguments 
that the function to be minimized is very natural. 
The relationship with some notions about functions of bounded variation developed by J.J.~Moreau 
is emphasized. 
\end{abstract} 

\section{Introduction} 
Let $\mathbf P$ a probability law on $\mathbb R$. 
Assume that $\mathbf P$ is of order $1$ (which writes $\int |x| \, d\mathbf P(x) < +\infty$). 
A real number $\bar m$ is a median if and only if it minimizes (references to 
exercises in some textbooks will be given in Section~\ref{Mq}) the function of $m$, 
$\int_{\mathbb R} |x - m| \, d\mathbf P(x)$. 
Without the order $1$ hypothesis, the medians minimize 

\[ 
m \mapsto \Phi(m) := \int_{\mathbb R} \bigl(|x - m| - |x|\bigr) \, d\mathbf P(x) \,. 
\] 
(This is not exactly, when $\tau = 1/2$ --- the value of $\tau$ corresponding to medians~---, 
the function defined in \eqref{Phi}, because of a factor $2$.) 
This extends to quantiles. 
This paper gives ``Convex Analysis'' proofs of these results. 
With our process the function $\Phi$ appears naturally. 
We emphasize the links with some notions about functions of bounded 
variation developed by J.J.~Moreau. 

I began this paper being unaware of Koltchinskii \cite{Kol}. 
This 1997 work is more devoted to the multivariate case. 
See some other comments in Section~\ref{convexview}. 
\medskip

\section{Definitions and notations}\label{DN} 
Let $\mathbf P$ be a probability law on $\mathbb R$. 
If necessary $X$ will denote a random variable obeying the law $\mathbf P$. 
The law $\mathbf P$ \emph{is of order} $p$ ($p \in \{1,2\}$) if 
$X$ is of order $p$, which writes $X \in L^p$. 
If no integrability condition on $X$ is assumed, specialists write $X \in L^0$ (order $0$). 
When $\mathbf P$ is of order $1$ the mean $\mathbb E(X)$ does exist. 
When $\mathbf P$ of order $0$ there exists a non-empty compact interval of 
\emph{medians} (the definition of a median is included in that of a quantile below). 
Let $F$ be the right-continuous distribution function: 
\[ 
F(x) := \mathbf P(\left]-\infty,x\right]) = \mathbf P({X \leq x}) \,. 
\] 
For $\tau \in \ ]0,1[$ (usually in Statistics $0.95$ or $0.99$ etc.) 
a real number $\bar q$ is a $\tau$-\emph{quantile} if 
\[ 
\mathbf P(X \leq \bar q) \geq \tau \quad \text{and} \quad 
\mathbf P(X \geq \bar q) \geq 1 - \tau 
\] 
or maybe more plainly 
\[ 
\mathbf P({X < \bar q}) \leq \tau \leq \mathbf P(X \leq \bar q) 
\] 
(note that $x \mapsto \mathbf P({X < x})$ is the left-continuous version $F^-$ of $F$). 
When $\tau = 1/2$ one recovers median(s). 
A geometrical definition of $\tau$-quantiles is the following: 
their set is the projection on $\mathbb R$ of the intersection of $\mathbb R \times \{\tau\}$ 
with the \emph{filled-in graph}\footnote{\ This graph is obtained from 
the graph of $F$ by adding vertical segments when there are gaps, thus obtaining 
an arcwise connected curve: 
\[ 
G = \{(x,y) \,;\, x \in \mathbb R \text{ and } y \in [F^-(x),F(x)]\} \,. 
\] 
This comes from Monteiro Marques and Moreau works on the sweeping process: 
\cite[p.147]{MM1}, \cite[p.15]{MM2}, some ideas coming back to \cite{Mor2}. 
In Section~\ref{convexview} the filled-in graph will be seen 
as a maximal monotone operator.} $G$ of $F$. 
That is, the set of $\tau$-quantiles is 
\[ 
\{\bar q \in \mathbb R \,;\, (\bar q,\tau) \in G \} \,. 
\] 
This geometrical definition operates when $\tau = 0$ or $1$. 
For example the set of $1$-quantiles is the closed, possibly empty, interval 
$\{q \in \mathbb R \,;\, F(q) = 1\}$ and the set of $0$-quantiles is the closed, possibly empty, interval 
$\{q \in \mathbb R \,;\, F^-(q) = 0\}$. 

A real number candidate to be the mean or a median will be denoted by $m$. 
A candidate to be a quantile will be denoted by $q$. 
\medskip

\section{Medians and quantiles as minimizers}\label{Mq} 
Well known\footnote{\ In French textbooks: \cite[Exercice E.11 p.71]{DCD1} 
(solution page 28 of \cite{DCD2}), and Th\'eor\`eme 10.1 page 93 in \cite{FF}. 
See also \cite[ex.1.8.2 p.51]{Fe1} and its solution in \cite{Fe2}.} 
is the result: \emph{If $\mathbf P$ is of order $1$, $\,\overline{\!m}$ is a median 
if and only if it minimizes on $\mathbb R$ the function} 
\[ 
m \mapsto \int_\mathbb R |x - m | \, d\mathbf P(x) \,. 
\] 
Surprising are: the assumption about order $1$, and the fact that medians depend only on the 
structure of ordered space of $\mathbb R$ and not on its metric (nor on its group structure, 
nor on Haar measure). 
An ``answer'', at least relatively to the ``order $1$'' assumption, is the following. 

Well known too\footnote{\ Exercises in French Universities and surely in 
many countries...} is the result: 
\emph{If $\mathbf P$ is any law, 
$x \mapsto |x - m | - |x|$ is $\mathbf P$-integrable 
(obviously $\bigl||x - m | - |x|\bigr| \leq |m|$) and $\,\overline{\!m}$ is a median 
if and only if it minimizes on $\mathbb R$ the function} 
\[ 
m \mapsto \int_\mathbb R \bigl(|x - m | - |x|\bigr) \, d\mathbf P(x) \,. 
\] 

The notion of median extends to $\mathbb R^d$ and to Banach spaces: see \cite{Kem,Kol,MD} 
and there exist conditional medians \cite{V3}; 
Kemperman \cite{Kem} and Milasevic \& Ducharme \cite{MD} 
gave in the multivariate case a sufficient condition implying uniqueness of the median 
(see already in 1948 Haldane \cite{Ha}). 
These questions will not be considered here. 

The interest of medians comes from robustness, i.e., stability with respect to outliers values. 
Maybe the notions of means in metric spaces going back to Fr\'echet 
(see \cite{Fr} and many other papers by the same author) should be revisited. 
See the papers by Armatte \cite{A1,A2}, the first one containing more than six 
pages of Fr\'echet's references. 
Numerous authors studied random variables with values in a metric space: 
for instance \cite{BH,RF}. 
\medskip 

Now we turn to quantiles. 
T.S. Ferguson \cite[Exercise 1.8.3 p.51, solution in \cite{Fe2}]{Fe1} says that\footnote{\ 
This result is appreciated by specialists: \cite[first chapter]{Koe} and already \cite[p.38]{KB}.}: 
\emph{If $\tau \in \ ]0,1[$ and $\mathbf P$ is of order $1$, $\,\overline{\!q}$ is a $\tau$-quantile 
if and only if it minimizes on $\mathbb R$ the function} 
\[ 
q \mapsto \int \rho_\tau(x - q) \, d\mathbf P(x) \,. 
\] 
\emph{where} 
\begin{equation}\label{rho} 
\rho_\tau(x) 
= \frac{1}{2} \, |x| + \bigl(\tau - \frac{1}{2}\bigr) \, x 
= \begin{cases} 
(\tau - 1) \, x &\text{if } x \leq 0, \\ 
\tau \, x &\text{if } x \geq 0. 
\end{cases} 
\end{equation} 
My purpose here is to give a \emph{Convex Analysis} proof of the extension to any 
probability law and to follow a natural way. 
Surely Koltchinskii \cite{Kol} contains the statement but his framework is multivariate. 
\medskip

\section{A Convex Analysis point of view}\label{convexview} 
The filled-in graph $G$ defined in Section~\ref{DN} is a subset of $\mathbb R^2$ 
which is a \emph{maximal monotone operator} \cite[Chapter~12]{RW} 
(and therefore a maximal cyclically monotone one because of dimension $1$, 
cf.\ \cite[12.6 pp.547--548]{RW}). 
Let $\mathbf F$ be the primitive (antiderivative) of $F$ null at $0$: 
\[ 
\mathbf F(x) = 
\begin{cases} 
\displaystyle \int_0^x F(u) \, du &\text{if } x \geq 0, \\ 
\displaystyle -\int_x^0 F(u) \, du &\text{if } x < 0 \,. 
\end{cases} 
\] 
Several textbooks treats convexity. 
One of the most fundamental is \cite{Mor1}. 
Since $F$ is nondecreasing $\mathbf F$ is convex. 
The \emph{sub-derivative} of $\mathbf F$ at $x$ is 
\[ 
\partial \mathbf F(x) 
= \{\ell \in \mathbb R \,;\, \forall h\in \mathbb R,\ \ell \, h \leq \mathbf F(x + h) \} \,. 
\]

\begin{lemma}\label{monotone} 
The graph of the multifunction $x \mapsto \partial \mathbf F(x)$ 
is nothing else but the filled-in graph $G$ of the graph of $F$: 
\[ 
G = \{(x,\ell) \in \mathbb R^2 \,;\, \ell \in \partial \mathbf F(x) \} \,. 
\] 
\end{lemma} 

\noindent 
{\sc Proof.} 
Let $\mathbf F'(x;w)$ denote the directional derivative of $\mathbf F$ at $x$ in the direction $w$. 
There holds 
\[ 
\partial \mathbf F(x) = [-\mathbf F'(x;-1),\mathbf F'(x;1)] \,. 
\] 
Then, since $F$ is right continuous, 
\begin{align*} 
\mathbf F'(x;1) &= \lim_{h \searrow 0} \frac{\mathbf F(x + h) - \mathbf F(x)}{h} \\ 
&= \lim_{h \searrow 0} \frac{1}{h} \, \int_x^{x+h} F(u) \, du \\ 
&= F(x) 
\end{align*} 
and, since $\mathbf F$ is as well the primitive of $F^-$, 
\begin{align*} 
\mathbf F'(x;-1) &= \lim_{h \searrow 0} \frac{\mathbf F(x - h) - \mathbf F(x)}{h} \\ 
&= -\lim_{h \searrow 0} \frac{1}{h} \, \int_{x-h}^x F^-(u) \, du \\ 
&= -F^-(x) \,. 
\end{align*} 
Therefore 
\[ 
\partial \mathbf F(x) = [F^-(x),F(x)] \,. \quad \Box 
\] 
\medskip 

In such a situation a natural question is: what does give the minimisation of $\mathbf F$? 
Obviously the infimum may be $-\infty$: this is why we cannot consider 
$\displaystyle \int_{-\infty}^x F(u) \, du$. 
But with a slope $\tau$ belonging to $]0,1[$ the function 
\[ 
q \mapsto \mathbf F(q) - \tau \, q 
\] 
does achieves minimum(s) and $\bar q$ is a minimum if and only if 
$0 \in -\tau + \partial \mathbf F(\bar q)$. 
This is equivalent to $(\bar q,\tau) \in G$, that is, $\bar q$ is a $\tau$-quantile. 
\medskip 

\noindent 
{\bf Scholium} \emph{Let $\mathbf F$ denote a primitive of $F$. 
For $\tau \in \left]0,1\right[$,} 
\begin{equation}\label{infimum} 
\bar q\, \text{ is a $\tau$-quantile} 
\Longleftrightarrow \bar q \text{ minimizes } q \mapsto \mathbf F(q) - \tau \, q \,. 
\end{equation}

Next Lemma is an integration by parts result. 
We follow a proof by Schilling but we could deduce the statement from 
results by Rockafellar \cite[Prop.1 pp.161--162]{Ro} or Moreau \cite[Section~11]{Mor3}. 
This will be detailed elsewhere.

\begin{lemma}\label{quantileplus} 
Let $a < b$ in $\mathbb R$. 
Then (note the half-open interval $]a,b]$) 
\begin{equation}\label{JJ} 
\int_{]a,b]} x \, d\mathbf P(x) = b \, F(b) - a \, F(a) - \int_a^b F(x) \, dx \,. 
\end{equation} 
\end{lemma} 

\noindent 
{\sc Proof} (from Schilling \cite[Exercise~13.13 p.133]{S1} 
and its solution on the Net \cite[pp.10--12]{S2}). 
For the product of the Lebesgue measure (we denote it by $dx$) and of $\mathbf P$ 
\begin{equation}\label{S1} 
(dx \otimes d\mathbf P)(]a,b]^2) = (b - a) \, [F(b) - F(a)] \,. 
\end{equation} 
And moreover 
\begin{align}\label{S2} 
(dx &\otimes d\mathbf P)(]a,b]^2) = \nonumber \\ 
&= \iint \mathbf 1_{]a,b]}(x) \, \mathbf 1_{]x,b]}(y) \, dx \, d\mathbf P(y) 
+ \iint \mathbf 1_{]a,b]}(x) \, \mathbf 1_{]a,x]}(y) \, dx \, d\mathbf P(y) \nonumber \\ 
&= \iint \mathbf 1_{]a,b]}(x) \, \mathbf 1_{]x,b]}(y) \, dx \, d\mathbf P(y) 
\overset{\text{Tonelli}}{+} \iint \mathbf 1_{]a,b]}(y) \, \mathbf 1_{[y,b]}(x) \, dx \, d\mathbf P(y) \nonumber \\ 
&= \int_{]a,b]} [F(b) - F(x)] \, dx + \int_{]a,b]} (b - y) \, d\mathbf P(y) \nonumber \\ 
&= (b - a) \, F(b) - \int_a^b F(x) \, dx + b \, [F(b) - F(a)] - \int_{]a,b]} y \, d\mathbf P(y) \,. 
\end{align} 
Comparing \eqref{S1} and \eqref{S2} one gets \eqref{JJ}. $\Box$

\medskip 
Now let\footnote{\ Koltchinskii \cite{Kol} considers the same functional. 
See specially two lines after his formula~(1.1) page~436 where he sets $f_{P,t}(s) := f_P(s) - s \, t$: 
there $s$ is my variable $q$ and $t$ my threshold $\tau$. 
Comparisons are not easy. 
He also uses (on pp.439--440) a multivalued integration result \cite[Th.8.3.4]{IT}. 
I got similar results earlier in \cite{V1,V2} but Ioffe and Tihomirov published also 
some papers on these questions before their book.} 
($\rho_\tau$ is defined in \eqref{rho} above) 
\begin{equation}\label{Phi} 
\Phi (q) := \int_{\mathbb R} \bigl(\rho_\tau(x - q) - \rho_\tau(x)\bigr) \, d\mathbf P(x) 
\end{equation} 
(note that the function $x \mapsto \rho_\tau(x - q) - \rho_\tau(x)$ is bounded on $\mathbb R$).

\begin{theo}\label{primitive} 
Let $\tau \in \left]0,1\right[$ and $\mathbf P$ be a law on $\mathbb R$. 
There holds 
\begin{equation*}\label{prim} 
\Phi(q) = 
\begin{cases} 
-\tau \, q + \int_{]0,q]} F(x) \, dx &\text{\rm if } q \geq 0, \\ 
-\tau \, q - \int_{]q,0]} F(x) \, dx &\text{\rm if } q < 0, 
\end{cases} 
\end{equation*} 
that is, $\Phi(q) = \mathbf F(q) -\tau \, q$. 
The function $\Phi$ is convex and inf-compact. 
The value $\bar q$ is a $\tau$-quantile if and only if it minimizes the function $\Phi$. 
\end{theo} 

\noindent 
{\sc Comments.} 
Part 3) below is a bit technical. 
Surely careless calculus would be quicker: the formal calculus 
\begin{align*} 
\Phi'(q) &= \int_\mathbb R \frac{d}{dq} \bigl(\rho_\tau(x - q) - \rho_\tau(x)\bigr) \, d\mathbf P(x) \\ 
 &= \int_{-\infty}^q (1 - \tau) \, d\mathbf P + \int_q^{+\infty} (-\tau) \, d\mathbf P \\ 
&= (1 - \tau) \, F(q) - \tau \, \bigl(1 - F(q)\bigr) \\ 
&= F(q)- \tau 
\end{align*} 
leads to 
\[ 
\Phi'(q) = 0 \Longleftrightarrow F(q) = \tau \ \text{!} 
\] 

\noindent 
{\sc Proof.} 
1) Firstly $x \mapsto \rho_\tau(x - q) - \rho_\tau(x)$ is less than 
$\tau \, |q|$ if $\tau \geq 1/2$ or than $(1 - \tau) \, |q|$ if $\tau \leq 1/2$. 
Or more roughly 
$\forall x$, $|\rho_\tau(x - q) - \rho_\tau(x)| \leq |q|$. 
So the integral is well defined. 
\medskip 

2) The function $\rho_\tau$ defined by \eqref{rho} being convex, 
the convexity of $q\mapsto \rho_\tau(x - q) - \rho_\tau(x)$ 
is trivial, hence the convexity of $\Phi$. 
\medskip 

3) We now prove $\Phi(q) \rightarrow +\infty$ as $|q| \rightarrow +\infty$. 

For $q \geq 0$, 
\[ 
\rho_\tau(x - q) - \rho_\tau(x) = 
\begin{cases} 
(1 - \tau) \, q &\text{if } x \in \left]-\infty,0\right], \\ 
(1 - \tau) \, q - x &\text{if } x \in \left]0,q\right], \\ 
-\tau \, q &\text{if } x \in \left]q,+\infty\right[, 
\end{cases} 
\] 
and, for $q < 0$, 
\[ 
\rho_\tau(x - q) - \rho_\tau(x) = 
\begin{cases} 
(1 - \tau) \, q &\text{if } x \in \left]-\infty,q\right], \\ 
x - \tau \, q &\text{if } x \in \left]q,0\right], \\ 
-\tau \, q &\text{if } x \in \left]0,+\infty\right[. 
\end{cases} 
\] 
Hence, if $q \geq 0$ (we will apply formula \eqref{JJ} between lines 2 and 3), 
\begin{align} 
\Phi(q) &= (1 - \tau) \, q \, \mathbf P(\left]-\infty,0\right]) + 
\int_{]0,q]} \bigl((1 - \tau) \, q - x\bigr) \, d\mathbf P(x) -\tau \, q \, \mathbf P(]q,+\infty[) \nonumber \\ 
&= (1 - \tau) \, q \, \mathbf P(\left]-\infty,q\right]) - \int_{]0,q]} x \, d\mathbf P(x) 
-\tau \, q \, \mathbf P(]q,+\infty[) \nonumber \\ 
&= (1 - \tau) \, q \, \mathbf P(\left]-\infty,q\right]) - \left[x \, F(x)\right]_0^q 
+ \int_{]0,q]} x \, F(x) \, dx 
-\tau \, q \, \mathbf P(]q,+\infty[) \nonumber \\ 
&= (1 - \tau) \, q \, F(q) - q \, F(q) + \int_{]0,q]} x \, F(x) \, dx 
-\tau \, q \, (1 - F(q)) \nonumber \\ 
&= -\tau \, q + \int_{]0,q]} F(x) \, dx \,. \label{Phi1} 
\end{align} 
When $q \rightarrow +\infty$ the mean value of $F$ over $]0,q]$ ultimately exceeds any 
value in $]\tau,1[$, so $\Phi(q) \rightarrow +\infty$. 
(More plainly for $q >\tau$, 
$\Phi(q) = \int_{]0,\tau]}\bigl(F(x)-\tau\bigr) \, dx + \int_{]\tau,q]}\bigl(F(x)-\tau\bigr) \, dx$ 
and $F(x) - \tau \rightarrow 1 - \tau >0$ as $x \rightarrow \infty$.) 
\medskip 

On the other hand, if $q < 0$, (we again apply \eqref{JJ} between lines 2 and 3), 
\begin{align} 
\Phi(q) &= (1 - \tau) \, q \, \mathbf P(\left]-\infty,q\right]) + 
\int_{]q,0]} (x - \tau \, q) \, d\mathbf P(x) -\tau \, q \, \mathbf P(]0,+\infty[) \nonumber \\ 
&= (1 - \tau) \, q \, \mathbf P(\left]-\infty,q\right]) + \int_{]q,0]} x \, d\mathbf P(x) 
-\tau \, q \, \mathbf P(]q,+\infty[) \nonumber \\ 
&= (1 - \tau) \, q \, \mathbf P(\left]-\infty,q\right]) + \left[x \, F(x\right]_q^0 
- \int_{]q,0]} x \, F(x) \, dx - \tau \, q \, \mathbf P(]q,+\infty[) \nonumber \\ 
&= (1 - \tau) \, q \, F(q) -q \, F(q) 
- \int_{]q,0]} x \, F(x) \, dx - \tau \, q \, (1 - F(q)) \nonumber \\ 
&= -\tau \, q - \int_{]q,0]} F(x) \, dx \,. \label{Phi2} 
\end{align} 
When $q \rightarrow -\infty$ the mean value of $F$ over $]q,0]$ ultimately passes under any 
value in $]0,\tau[$, so $\Phi(q) \rightarrow +\infty$. 
\medskip 

Thus the finite valued convex function $\Phi$ is inf-compact, so it achieves its 
infimum over a non-empty compact interval. 
By \eqref{Phi1} and \eqref{Phi2} 
\[ 
\Phi(q) = \mathbf F(q) -\tau \, q \,. 
\] 
From \eqref{infimum} the minimizers are the $\tau$-quantiles. 
$\Box$

\begin{cor}\label{quantile} 
Let $\tau \in \left]0,1\right[$ and $\mathbf P$ be a first order law on $\mathbb R$ then 
\[ 
\bar q \text{ is a $\tau$-quantile} 
\Longleftrightarrow 
\bar q \text{ minimizes } 
q \mapsto \int_{\mathbb R} \rho_\tau(x - q) \, d\mathbf P(x) \,. 
\] 
\end{cor} 
 
\noindent 
{\sc Proof.} 
The term 
\[ 
\int_{\mathbb R} \rho_\tau(x) \, d\mathbf P(x)
\] 
is finite and does not depend on $q$. 
Hence one can add it to the right-hand side of \eqref{Phi}.
$\Box$ 
\medskip 

\section{A direct proof} 

Classically the convex function $\Phi$ on $\mathbb R$ achieves a minimum at $\bar q$ 
if and only if the left and right derivatives of $\Phi$ at $\bar q$ are respectively 
$\leq 0$ and $\geq 0$.

\begin{theo}\label{quantileplus} 
Let $\tau \in \left]0,1\right[$ and $\mathbf P$ be a law on $\mathbb R$. 
The right derivative of $\Phi$ at $\bar q$ is $\geq 0$ if and only if 
$\mathbf P(\left]-\infty,\bar q\,\right])$ is $\geq \tau$. 
The derivative of $\Phi$ at $\bar q$ in direction $-1$ is $\geq 0$ if and only if 
$\mathbf P(\left]-\infty,\bar q\right[)$ is $\leq \tau$. 
\end{theo} 

\noindent 
{\sc Proof.} 
We have to reformulate the inequalities 
\[ 
\Phi'(\bar q;1) = \lim_{h \searrow 0} \frac{\Phi(\bar q + h) - \Phi(\bar q)}{h} \geq 0 
\] 
and 
\[ 
\Phi'(\bar q;-1) = \lim_{h \searrow 0} \frac{\Phi(\bar q - h) - \Phi(\bar q)}{h} \geq 0 
\] 
the second one expressing the positiveness of the directional derivative at $\bar q$ 
in the direction $-1$ (another way to say that the left derivative is $\leq 0$). 
\medskip 

a) Firstly let us consider the right derivative at $\bar q$. 
One has 
\[ 
\frac{\Phi(\bar q + h) - \Phi(\bar q)}{h} = 
\frac{1}{h}\int_{\mathbb R} \bigl(\rho_\tau(x - \bar q - h) - \rho_\tau(x - \bar q)\bigr) \, d\mathbf P(x) \,. 
\] 
The integral splits into 
\[ 
\int_{\left]-\infty,\bar q\right]} + \int_{]\bar q,\bar q + h[} + \int_{[\bar q + h,+\infty[} 
\] 
where the (non written) integrands are respectively $(1 - \tau)\,h$, 
a function bounded above by $h$ (on a \emph{vanishing domain}\footnote{\ The open intervals 
$]\bar q,\bar q + h[$ decrease when $h \searrow 0$ and have empty intersection. 
In order to apply convergence theorems of Integration Theory one should consider 
--- and this is sufficient --- a sequence $(h_n)_{n\in\mathbb N}$ satisfying $h_n \searrow 0$.}) 
and $-\tau \, h$ (on a domain which converges to $]\bar q,+\infty[$). 
Hence
\begin{align*} 
\lim_{h \searrow 0} \frac{\Phi(\bar q + h) - \Phi(\bar q)}{h} 
&= (1-\tau) \, \mathbf P(\left]-\infty,\bar q\right]) + (-\tau)\, \mathbf P(]\bar q,+\infty[) \\ 
&= \mathbf P(\left]-\infty,\bar q\right]) - \tau \,. 
\end{align*} 
Thus the right derivative is $\geq 0$ if and only if 
$\mathbf P(\left]-\infty,\bar q\right])$ is $\geq \tau$. 
\medskip 

b) We turn now\footnote{\ Minus signs are always more perilous.} 
to the directional derivative in the direction $-1$. 
One has 
\[ 
\frac{\Phi(\bar q - h) - \Phi(\bar q)}{h} = 
\frac{1}{h}\int_{\mathbb R} \bigl(\rho_\tau(x - \bar q + h) - \rho_\tau(x - \bar q)\bigr) \, d\mathbf P(x) \,. 
\] 
The integral splits into 
\[ 
\int_{\left]-\infty,\bar q - h\right]} + \int_{]\bar q - h,\bar q[} + \int_{[\bar q,+\infty[} 
\] 
where the integrands are respectively $-(1 - \tau)\,h$ 
(over a domain which converges to $\left]-\infty,\bar q\right[$), 
a function bounded above by $h$ (the domain being vanishing) and $\tau \, h$. 
Whence 
\begin{align*} 
\lim_{h \searrow 0} \frac{\Phi(\bar q - h) - \Phi(\bar q)}{h} 
&= -(1-\tau) \, \mathbf P(\left]-\infty,\bar q\right[) + \tau\, \mathbf P([\bar q,+\infty[) \\ 
&= -\mathbf P(\left]-\infty,\bar q\right[) + \tau \,. 
\end{align*} 
Thus the directional derivative in the direction $-1$ is $\geq 0$ if and only if 
$\mathbf P(\left]-\infty,\bar q\right[)$ is $\leq \tau$. 
$\Box$ 
\medskip

\noindent 
{\sc Acknowledgments.} Thanks to Manuel Monteiro Marques, Paul Raynaud de Fitte 
and Lionel Thibault for their comments. 
They are not responsible of possible weaknesses.

\providecommand{\bysame}{\leavevmode\hbox to3em{\hrulefill}\thinspace} 
\providecommand{\MR}{\relax\ifhmode\unskip\space\fi MR } 
\providecommand{\MRhref}[2]{%
 \href{http://www.ams.org/mathscinet-getitem?mr=#1}{#2} 
}
\providecommand{\href}[2]{#2}

\end{document}